\documentclass[11pt, reqno]{amsart}
\usepackage{a4wide,amsfonts,amsmath,amssymb,amsthm,epsfig,cite,graphicx,hyperref,
	color, esint,fancyhdr, enumerate, latexsym, amsrefs}

\numberwithin{equation}{section}

\usepackage{tikz-cd}

\newtheorem{thm}{Theorem}[section]

\theoremstyle{definition}

\hypersetup{
	colorlinks,
	citecolor=blue,
	filecolor=black,
	linkcolor=red,
	urlcolor=black
}

\newcommand{\D}{\mathbb{D}}
\newcommand{\C}{\mathbb{C}} 
\newcommand{\R}{\mathbb{R}}

\begin{document}
\title{A submultiplicative property of the Carath\'eodory metric on planar domains}
	
\author{Amar Deep Sarkar and Kaushal Verma}
	
\address{ADS: Department of Mathematics, Indian Institute of Science, Bangalore 560012, India}
\email{amarsarkar@iisc.ac.in}

\address{KV: Department of Mathematics, Indian Institute of Science, Bangalore 560 012, India}
\email{kverma@iisc.ac.in}
	
	

\begin{abstract}
Given a pair of smoothly bounded domains $D_1, D_2 \subset \mathbb C$, the purpose of this note is to obtain an inequality that relates the Carath\'{e}odory metrics on $D_1, D_2, D_1 \cap D_2$ and $D_1 \cup D_2$.
\end{abstract}
	
\maketitle
	
	
\section{Introduction}
\noindent Let $ \lambda_D(z) \vert dz \vert $ denote the Poincar\'e metric on a hyperbolic domain $ D \subset \C $. To quickly recall the construction of this metric, note that there exists a holomorphic covering from the unit disc $ \D $ to $D$,
\[
\pi : \D \rightarrow D
\] 	
whose deck transformations form a Fuchsian group $ G $ that acts on $ \D $. The Poincar\'e metric on $ \D $, 
\[
\lambda_{\D}(z)\vert dz \vert = \frac{\vert dz \vert}{1 - \vert z \vert^2}
\]
is invariant under ${\rm Aut}(\D) $ (hence in particular $ G $) and therefore, for $z \in D$, the prescription 
\begin{equation}\label{Eq:PullBackPoincare}
\lambda_D(\pi(z)) \vert \pi^{\prime}(z) \vert = \lambda_{\D}(z)  
\end{equation}
defines the Poincar\'e metric $ \lambda_D $ on $ D $ in an unambiguous manner. For an arbitrary $ z \in D $, we may choose the covering projection so that 
$ \pi(0) = z $ and hence (\ref{Eq:PullBackPoincare}) implies that
\[
\lambda_D(z) = \vert \pi^{\prime}(0) \vert^{-1}.
\] 
Solynin \cite{SolyninMultiplicativity}, \cite{Solynin2} proved the following remarkable relation between the Poincar\'e metrics on a pair of hyperbolic domains and those on their union and intersection: 

\medskip

\noindent {\it Let $ D_1,  D_2 \subset \C $ be domains such that $ D_1 \cap D_2 \neq \emptyset $. Suppose that $ D_1 \cup D_2 $ is hyperbolic. Then 
\[
\lambda_{D_1 \cap D_2}(z) \lambda_{D_1 \cup D_2}(z)  \leq \lambda_{D_1}(z) \lambda_{D_2}(z)
\]
for all $ z \in D_1 \cap D_2 $. If equality holds at one point $ z_0 \in D_1 \cap D_2$, then $ D_1 \subset D_2 $ or $ D_2 \subset D_1 $ and in this case, equality holds for all points $ z \in D_1 \cap D_2 $.}

\medskip

A direct proof of this was given by Kraus-Roth \cite{OliverRothMultiplicativity} that relied on a computation reminiscent of the classical Ahlfors lemma and the fact that 
\begin{equation}
	\Delta \log \lambda_D(z) = 4 \lambda_D^2(z),
\end{equation}
which precisely means $ \lambda_D(z) \vert dz \vert $ has constant curvature $ -4 $ on $ D $. Solynin's result follows from a comparison result for solutions to 
non-linear elliptic PDE's of the form
\begin{equation}\label{Eq:EllipticPDE}
\Delta u - \mu(u) - f = 0,
\end{equation}
where $ \mu : \R \rightarrow \R $ and $ f : D \rightarrow \R $ are suitable continuous non-negative functions that satisfy an additional convexity condition. Clearly, this reduces to the curvature equation (\ref{Eq:EllipticPDE}) by writing $ u = \log \lambda_D $ and letting $ \mu(x) = 4 e^{2x} $ and $ f \equiv 0 $.

\medskip

The purpose of this note is to prove an analogue of Solynin's theorem for the Carath\'{e}odory metric on planar domains. Let $ D \subset \C $ be a domain that admits at least one non-constant bounded holomorphic function. Recall that for 
$ z \in D $, the Carath\'eodory metric $ c_D(z) \vert dz \vert $ is defined by 
\[
c_D(z) \vert dz \vert = \sup\{\vert f^{\prime}(z) \vert : f : D \rightarrow \D \,\, \text{holomorphic and} \, \, f(z) = 0\}.
\]
This is a distance-decreasing (and hence conformal) metric in the sense that if $h : U \rightarrow V$ is a holomorphic map between a pair of planar domains $U, V$, then 
\[
c_V(h(z)) \lvert h'(z) \rvert \le c_U(z), \;\; z \in U.
\]
If $U \subset V$, applying this to the inclusion $i : U \rightarrow V$ shows that the Carath\'{e}odory metric is monotonic as a function of the domain, i.e., $c_V(z) \le c_U(z)$ for $z \in U$. Furthermore, for each $\zeta \in D$, there is a unique holomorphic map $f_{\zeta} : D \rightarrow \D$ that realizes the supremum in the definition of $c_D(z)$ -- this is the Ahlfors map. 

\section{Statement and proof of the main result}
\begin{thm}
	Let $ D_1, D_2 \subset \C$ be smoothly bounded domains with $ D_1 \cap D_2 \neq \emptyset $. Then there exists a constant $ C = C(D_1, D_2) > 0 $ such that 
	\[
	c_{D_1 \cap D_2}(z) \cdot c_{D_1 \cup D_2}(z) \leq C\, c_{D_1}(z) \cdot c_{D_2}(z)
	\]
	for all $z \in D_1 \cap D_2$.
\end{thm} 
Note that $D_1 \cap D_2$ may possibly have several components. The notation $c_{D_1 \cap D_2}(z)$ refers to the Carath\'{e}odory metric on that component of $D_1 \cap D_2$ which contains a given $z \in D_1 \cap D_2$. Furthermore, the constant $C > 0$ is independent of which component of $D_1 \cap D_2$ is being considered and only depends on $D_1$ and $D_2$. The proof, which is inspired by \cite{OliverRothMultiplicativity}, uses the following known supplementary properties of the Carath\'{e}odory metric:

\medskip

First, $c_D(z)$ is continuous and $\log c_D(z)$ is subharmonic on $D$ -- see \cite{AhlforsAndBeurling} for instance. The possibility that $\log c_D(z) \equiv -\infty$ can be ruled out for bounded domain since for every $\zeta \in D$, the affine map $f(z) = (z - \zeta)/a$ vanishes at $\zeta$ and maps $D$ into the unit disc $\mathbb D$ for every positive $a$ that is bigger than the diameter of $D$. The fact that $f'(z) \equiv 1/a > 0$ shows that $c_D(\zeta) > 0$ and hence 
$\log c_D(\zeta) > -\infty$.  

\medskip

Second, Suita \cite{SuitaI} showed that $c_D(z)$ is real analytic in fact and hence we may speak of its curvature 
\[
\kappa_D(z) = -c_D^{-2}(z) \Delta \log c_D(z)
\]  
in the usual sense. The subharmonicity of $\log c_D(z)$ already implies that $\kappa_D \le 0$ everywhere on $D$, but by using the method of supporting metrics, Suita \cite{SuitaI} was also able to prove a much stronger inequality namely, 
$\kappa_D \le -4$ on $D$. Following this line of inquiry further, Suita \cite{SuitaII} (see \cite{BurbeaPaper} as well) showed that if the boundary of $D$ consists of finitely many Jordan curves, the assumption that $\kappa_D(\zeta) = -4$ for some point $\zeta \in D$, implies that $D$ is conformally equivalent to $\mathbb D$. 

\medskip

Finally, it is known that this metric admits a localization near $C^{\infty}$-smooth boundary points -- see for example \cite{InvariantMetricJarnicki} which contains a proof for the case of strongly pseudoconvex points that works verbatim in the planar case too. That is, if $p \in \partial D$ is a $C^{\infty}$-smooth boundary point of a bounded domain $D \subset C$, then for a small enough neighbourhood $U$ of $p$ in $\mathbb C$,
\[
\lim_{z \rightarrow p} \frac{c_{U \cap D}(z)}{c_D(z)} = 1.
\] 
Using this, it was shown in \cite{SV} that the curvature $\kappa_D(z)$ of the 
Carath\'{e}odory metric approaches $-4$ near each $C^{\infty}$-smooth boundary point of a bounded domain $D \subset \mathbb C$. The point here being that as one moves nearer to such a point, the metric begins to look more and more like the Carath\'{e}odory metric on $\mathbb D$ -- the use of the scaling principle makes all this precise. It follows on any smoothly bounded planar domain $D$, 
$ \kappa_D \approx -4 $ for points close to the boundary and for those that are at a fixed positive distance away from it, there is a large negative lower bound for $\kappa_D$ as a result of its continuity. Hence for every such $ D \subset \C$, there is a constant $ C = C(D) > 0$ such that 
\[
-C \leq \kappa_D(z) \leq -4
\]
for all $ z \in D $. Another consequence of the localization principle is that $c_D(z) \rightarrow +\infty$ as $z$ approaches a smooth boundary point. Indeed, near such a point, $c_D(z) \approx c_{U \cap D}(z)$ and the latter is the same as the Carath\'{e}odory metric on the disc $\mathbb D$ (since $U \cap D$ can be taken to be simply connected) which blows up near every point on $\partial \mathbb D$.

\medskip

To prove the theorem, let $ \kappa_1 $ and $ \kappa_2 $ be the curvatures of the Carath\'eodory metric $ c_{D_1}(z) \lvert dz \rvert $ and $ c_{D_2}(z) \lvert dz \rvert $ on $ D_1 $ and $ D_2 $ respectively.

\medskip

Consider the metric
\[
c(z) \lvert dz \rvert = \frac{c_{D_1}(z)\cdot c_{D_2}(z)}{c_{D_1 \cup D_2}(z)} \lvert dz \rvert
\]
on the possibly disconnected open set $ D_1 \cap D_2 $. What follows applies to each component of $D_1 \cap D_2$ without any regard to its analytic or topological properties and hence we will continue to write $c(z) \lvert dz \rvert$ to denote this metric on any given component. Its curvature is
\begin{equation}
\begin{split}
\kappa(z) &= -c^{-2}(z) \Delta \log c(z)\\
 &= -c^{-2}(z)\left( \Delta \log c_{D_1}(z) + \Delta \log c_{D_2}(z) - \Delta \log c_{D_1 \cup D_2}(z) \right)\\
 &= I_1 + I_2 + I_{12},
\end{split}
\end{equation}
where
\[
I_1 = -c^{-2}(z) \Delta \log c_{D_1}(z),
\]
\[
I_2 = -c^{-2}(z) \Delta \log c_{D_2}(z) 
\]
and
\[
I_{12} = c^{-2}(z)\Delta \log c_{D_1 \cup D_2}(z).
\]
Note that $ I_{12} \geq 0$ since $ \log c_{D_1 \cup D_2}(z) $ is subharmonic and hence
\[
\kappa(z) \geq I_1 + I_2.
\]
To analyze each of these terms, note that $ c_{D_1} \geq c_{D_1 \cup D_2} $ and $ c_{D_2} \geq c_{D_1 \cup D_2} $. Combining this with the fact that the curvature of the Carath\'eodory metric is negative everywhere, it follows that 
\[
I_1 \geq \kappa_1(z)
\] 
and 
\[
I_2 \geq \kappa_2(z).
\] 
Hence there is a constant $ C = C(D_1, D_2) > 0 $ such that
\[
\kappa(z) \geq \kappa_1(z) + \kappa_2(z) > -C
\]
for all $ z \in D_1 \cap D_2 $.

\medskip

Let $U$ be a component of $D_1 \cap D_2$. Since $D_1, D_2$ have smooth boundaries, $U$ has finite connectivity, say $m \ge 1$ and is non-degenerate in the sense that the interior of its closure coincides with itself. In particular, its boundary cannot contain isolated points. Let $U_{\epsilon}$ be an $\epsilon$-thickening of $U$. For all sufficiently small $\epsilon > 0$, $U_{\epsilon}$ also has connectivity $m$ and is non-degenerate. Furthermore, $U_{\epsilon} \rightarrow U$ in the sense of Carath\'{e}odory as $\epsilon \rightarrow 0$.

\medskip

For a fixed $\epsilon > 0$, let $c_{\epsilon}(z) \vert dz \vert$ be the 
Carath\'{e}odory metric on $U_{\epsilon}$. Consider the function
\[
u(z) = \log\left(\frac{c_{\epsilon}(z)}{\sqrt{C\big/ 4 }\,c(z)}\right)
\]
on $U$. Since $ U$ is compactly contained in $ U_{\epsilon} $, 
$ c_{\epsilon} $ is bounded on $U$ and if $ \xi \in \partial U $, then as 
$z \rightarrow \xi$ within $U$,  
\[
\lim_{z \to \xi} \frac{c_{D_1}(z)\cdot c_{D_2}(z)}{c_{D_1 \cup D_2}(z)} \geq \lim_{z \to \xi}c_{D_2}(z) = + \infty,
\]
where the inequality follows from the monotonocity of the Carath\'{e}odory metric, i.e., $ c_{D_1}(z) \geq c_{D_1 \cup D_2}(z) $. The fact that $\xi$ is a smooth boundary point of $D_2$ implies that $c_{D_2}$ blows up near it and this means that  $ u(z) \to - \infty $ at $\partial U$. Therefore, $u$ attains a maximum at some point $z_0 \in U$. As a result, 
\begin{equation}
\begin{split}
0 \geq \Delta u(z_0) &= \Delta \log c_{\epsilon}(z_0) - \Delta \log c(z_0)\\
&\geq - \kappa_{{\epsilon}}(z_0)c_{\epsilon}^2(z_0) - C c^2(z_0)\\
&\geq 4 c^{2}_{\epsilon}(z_0) - C c^{2}(z_0),
\end{split}
\end{equation}
where $ \kappa_{\epsilon} $ is the curvature of $c_{\epsilon}(z) \vert dz \vert$.
It follows that
\[
u(z_0) = \log \left(\frac{c_{\epsilon}(z_0)}{\sqrt{C \big/4} \, c(z_0)}\right) \leq 0.
\]
For an arbitrary $ z \in U$, $ u(z) \leq u(z_0) \leq 0 $ and this gives 
\[
\log\left(\frac{c_{\epsilon}(z)}{\sqrt{C \big/4} \, c(z)}\right) \leq 0
\]
or
\[
c_{\epsilon}(z) \leq \sqrt{C \big/4}\, c(z)
\]
which is same as 
\[
c_{\epsilon}(z)\cdot c_{D_1 \cup D_2}(z) \leq \sqrt{C \big/4}\,c_{D_1}(z)\cdot c_{D_2}(z)
\]
and this holds for all $z \in U$. It remains to show that $c_{\epsilon} \rightarrow c_U$ pointwise on $U$ for then we can pass to the limit as $\epsilon \rightarrow 0$, keeping in mind that $C$ is independent of $\epsilon$, to get
\[
c_U(z)\cdot c_{D_1 \cup D_2}(z) \leq \sqrt{C \big/4}\,c_{D_1}(z)\cdot c_{D_2}(z)
\]
as claimed. Fix $p \in U$. To show that $c_{\epsilon}(p) \rightarrow c_U(p)$, it suffices to prove that 
$\vert f'_{\epsilon}(p) \vert \rightarrow \vert f'(p) \vert$ as $\epsilon \rightarrow 0$, where $f_{\epsilon} : U_{\epsilon} \rightarrow \mathbb D$ and
 $f : U \rightarrow \mathbb D$ are the respective Ahlfors maps for the domains $U_{\epsilon}$ and $U$ at $p$. What this is really asking for is the continuous dependence of the Ahflors maps on the domains. But this is addressed in \cite{Younsi} -- indeed, Theorem 3.2 therein can be applied as the domains $U_{\epsilon}, U$ have the same connectivity by construction and the 
 $U_{\epsilon}$'s decrease to $U$ as $\epsilon \rightarrow 0$. The nuance, in this theorem, about the base point being the point at infinity can be arranged by sending $p \mapsto \infty$ by the map $T(z) = 1/(z-p)$ and working with the domains $T(U_{\epsilon})$ and $T(U)$. This completes the proof.

\medskip

\noindent {\it Question:} Let $\Omega \subset \mathbb C$ be a domain on which the Carath\'{e}odory metric $c_{\Omega}(z) \vert dz \vert$ is not degenerate. Does its curvature $\kappa_{\Omega}$ admit a global lower bound?


\end{document}